\pdfoutput=1
\documentclass[review,12pt]{elsarticle}
\usepackage{amsmath}
\usepackage{amsfonts}
\usepackage{amssymb}
\usepackage{bbm}
\usepackage{bbold}
\usepackage{amsthm}
\usepackage{soul}
\usepackage{color}
\usepackage{cancel}

\usepackage{thmtools}
\usepackage{mathrsfs}

\theoremstyle{plain}
\newtheorem{theorem}{Theorem} 
\newtheorem{lemma}[theorem]{Lemma}

\newtheorem{proposition}[theorem]{Proposition}

\newtheorem{definition}[theorem]{Definition}

\newtheorem{example}[theorem]{Example}

\biboptions{authoryear}

\begin{document}

\begin{frontmatter}

\title{A technical note on  finite best-worst random utility model representations\\ Revised 16-09-2024}


\author{Hans Colonius}
\address{University of Oldenburg}

\ead[url]{https://uol.de/en/hans-colonius/}

\ead{hans.colonius@uni-oldenburg.de}


\begin{abstract}
This paper investigates the random utility representation of best-worst choice probabilities (picking the best and the worst alternative from an offered set). \cite{doignon2023} presented a complete characterization of the best-worst-choice polytope on four alternatives. Moreover, using polytope  methods he showed that the Block-Marschak inequalities for best-worst choices are not sufficient for a random utility representation of best-worst choices for sets of four or more alternatives. Following the approach of \cite{falmagne1978}, we construct a probability measure on the set of rankings for a set of four alternatives implying a random utility representation.

\end{abstract}


\end{frontmatter}

\section*{Introduction}
The purpose of this note is to construct a probability measure on the rankings (permutations) of a finite set $A$ that allows one to infer the existence of a random utility representation of best-worst choice probabilities for $A$. For context and related results, we refer to  \cite{marley2005}, \cite{colonius2021}, and \cite{doignon2023}.
\section{Some definitions and basic results}
For ease of reference, we keep the notation in the original paper of  \cite{falmagne1978}. For a finite set $A$, we write $|A|$ for the number of elements in $A$, $\mathscr{P}(A)$ for the power set of $A$. For  any finite  nonempty set $A$  let  $\Phi(A)$ be the set of all finite, nonempty subsets of $A$. Moreover, let  $\mathscr{P}(A,\geq c)$ (resp., $\mathscr{P}(A,\leq c)$) be the set of all subsets of $A$ containing at least (resp., at most) $ c$ elements. 
\begin{definition}\normalfont{
		Let $A$ be a nonempty set  of elements with at least 2 elements. For any $B\subseteq A$, with  $|B| \ge 2$, and any $a, b \in B$, $a\ne b$, 
		let $BW_B(a,b)\mapsto [0,1]$ denote  the probability that $a$ and $b$ are respectively chosen as best and worst elements in the subset $B$ of $A$. Let $\mathbb{BWP}= \{BW_B\,|\, B\subseteq A, BW_B(a,b)\mapsto [0,1],  a,b \in B, a\ne b\}$ be the collection of all those probability measures. Suppose that
		\[ \sum_{a,b \,\in B, \, a\ne b}  BW_B(a,b)=1 \hspace{1cm} (BW \in \mathbb{BWP}).\]
		Then $(A,\mathbb{BWP})$ is called a \emph{finite system of best-worst choice probabilities}, or more briefly, a \emph{system}. }
\end{definition}
The quantities $BW_B(a,b)$ will be referred to as \emph{best-worst choice probabilities} indicating that an alternative $a$ is judged as best and $b$ as worst, when subset $B$ is available. 
\subsection{Best-Worst Block-Marschak polynomials}
Suppose $(A,\mathbb{BWP})$ is a system. Consider the following expressions:
\begin{align*}
& BW_A(a,b),\\
& BW_{A\setminus\{c\}}(a,b)- BW_A(a,b),\\
& BW_{A\setminus\{c,d\}}(a,b)-[BW_{A\setminus\{c\}}(a,b)+BW_{A\setminus\{d\}}(a,b)]+BW_A(a,b),\\
&  BW_{A\setminus\{c,d,e\}}(a,b)-[BW_{A\setminus\{c,d\}}(a,b)+BW_{A\setminus\{c,e\}}(a,b)+BW_{A\setminus\{d,e\}}(a,b)]+\\
&\hspace{2cm}[BW_{A\setminus\{c\}}(a,b)+BW_{A\setminus\{d\}}(a,b)+BW_{A\setminus\{e\}}(a,b)]- BW_A(a,b),\\
& \text{etc.}
\end{align*}
%
We introduce a compact notation.
\begin{definition}
	For any $B\subset A$, $B\neq A$, $a,b \in A\setminus B, a\ne b,$ in a system $(A,\mathbb{BWP})$, we define
	\begin{equation}\label{BM}
	K_{ab,B} = \sum_{i=0}^{|B|} (-1)^i \sum_{C\in \mathscr{P}(B, |B|-i)} BW_{A\setminus C}(a,b).
	\end{equation}
	The $K_{ab,B}$ are called \emph{best-worst Block-Marschak polynomials of} $(A,\mathbb{BWP})$, or \emph{best-worst BM polynomials}, for short.
\end{definition}
Observe that 
\begin{align*}
K_{ab,\emptyset} &=\sum_{i=0}^0 (-1)^i \sum_{C\in \mathscr{P}(\emptyset, 0-i)} BW_{A\setminus C}(a,b)\\
& = BW_A(a,b);\\
K_{ab,\{c\}}    &=    BW_{A\setminus\{c\}}(a,b)- BW_A(a,b)  \\
& =    BW_{A\setminus\{c\}}(a,b)- K_{ab,\emptyset};\\
K_{ab,\{c,d\}} & =  BW_{A\setminus\{c,d\}}(a,b)-[BW_{A\setminus\{c\}}(a,b)+BW_{A\setminus\{d\}}(a,b)]+BW_A(a,b) \\
&= BW_{A\setminus\{c,d\}}(a,b)- K_{ab,\{c\}} - K_{ab,\{d\}} -K_{ab,\emptyset} .
\end{align*}
Similar computations show that 
\begin{align*}
K_{ab,\{c,d,e\}} &=  BW_{A\setminus\{c,d,e\}}(a,b)-  K_{ab,\{c,d\}} -  K_{ab,\{c,e\}} -  K_{ab,\{d,e\}} \\
&  \hspace{1cm}-K_{ab,\{c\}} - K_{ab,\{d\}}  - K_{ab,\{e\}}  -K_{ab,\emptyset}
\end{align*}
These examples suggest the following result.
\begin{theorem}\label{moebius}
	Let $(A,\mathbb{BWP})$ be  a system of best-worst choice probabilities. Then, for all $B\subset A, B\ne A$, and $a,b \in A\setminus B$,
	\begin{equation}
	BW_{A\setminus B}(a,b) = \sum_{C\in \mathscr{P}(B)} K_{ab, C}.
	\end{equation}
\end{theorem}

For proof see \cite{colonius2021} and the references therein.                                                                                                                                                  
%
%
\subsection*{Some more definitions}
%
Let $\ge$ be an arbitrarily chosen simple order on $A$. As usual, we write, for any $a,b \in A, a<b$ iff not $a\ge b$, and $a>b$ iff not $b \ge a$. For any $B\subset A$, we write $\Pi_B$ for the set of $|B|!$ permutations (rankings) on $B$. For $a_1, a_2,\ldots, a_n \in A$, we write $a_1a_2\ldots a_n$ for the ranking (defined by $>$) corresponding to permutation $\pi$ with $  \pi(a_1)>\pi(a_2)>\cdots  >\pi(a_n)$.

\begin{definition}
	Let  $B\subset A$ with $|A|-|B|\geq 2$; for $a, b \in A\setminus B$ define
				\[\mathrm{S}(a~ A\setminus B~b)=  \{\pi \in \Pi_A\,|\,    \text{ $a$ is ranked ``best''  and $b$ is ranked ``worst''  in $A\setminus B$}  \}.\]
\end{definition}
\begin{lemma}\label{lemma:disjoint}	
	For any distinct  $a,b \in A\setminus B$
	\begin{equation}
		S(a~A\setminus B~b) = \sum_{C\in \mathscr{P}(B)}\sum_{\pi\in \Pi_C} \sum_{(\pi_1 \pi_2)=\pi} S(\pi_1 a ~A \setminus C~ b \pi_2);
	\end{equation}
\end{lemma}
Here, notation  $S(\pi_1 a ~A \setminus C~ b \pi_2)$ means the set of all rankings of  $\Pi_A$ where $\pi=(\pi_1 \pi_2) \in \Pi_C$  is split into all possible pairs $\pi_1,\pi_2$; for example, ranking $cde\in \Pi_{\{c,d,e\}}$ is split  into \[ cde/\cdot, cd/e, c/de, \cdot/cde\] (using notation $\cdot /\cdot$) so that $cd/e$ corresponds to $S(cda~A\setminus C~be)$, $\cdot/cde$ to $S(a~A\setminus C~bcde)$, etc. For a proof of  Lemma \ref{lemma:disjoint} see \cite{colonius2021}.

\begin{example}
	Let $A=\{a,b,c,d\}$ and $B=\{c,d\}$. Then 
\begin{align*}
 S(a~A\setminus B~b) &=  	S(a~A~b) + S(ca~A~b) +  S(a~A~bc) + S(da~A~b) +  S(a~A~bd) \\
 +& S(cda~A~b)+S(dca~A~b) +S(a~A~bcd)+S(a~A~bdc)\\
 &+S(ca~A~bd)+S(da~A~bc)\\
 &= \{acdb,adcb\}+ \{cadb\}+\{adbc\}+\{dacb\}+\{acbd\}\\
 +&\{cdab\}+\{dcab\}+\{abcd\}+\{abdc\}+\{cabd\}+\{dabc\}.
\end{align*}
\end{example}

\section{Construction of  $\mathrm{P}[\,.\,]$}

The following lemma  is a variant of the well-known result by \cite{block1960} (see also \cite{marley2005}, section on random ranking models).
\begin{lemma}\label{lem:BM}
	$(A,\mathbb{BWP})$ is a best-worst random utility system if and only if there exists 	a probability measure $\mathrm{P}[\,.\,]$ on $\mathscr{P}(\Pi_A)$ satisfying
	\begin{equation}\label{extBM}
		BW_{B}(a,b) = \mathrm{P}[S(a~B~b)] 
	\end{equation}
	for all $a.b \in B\in \mathscr{P}(A,\geq 2)$.
\end{lemma}
In constructing a probability measure $\mathrm{P}[\,.\,]$ on $\mathscr{P}(\Pi_A)$ we must show, in addition to Equation~(\ref{extBM}), that $\mathrm{P}[\Pi_A]=1$. We first define $\mathrm{P}[\,.\,]$  on the points of $\Pi_A$ and then extend it in the usual way to a probability measure on $\mathscr{P}(\Pi_A)$.

 We define the probability of each point $\pi$ in $\Pi_A$ as follows: 
 \begin{definition}
 For $a, b \in A\setminus B$ 
 \begin{equation}\label{def}
	\mathrm{P}[\pi \in 	S(a~A\setminus B~b)] = \frac{K_{ab, B}}{|S(a~A\setminus B~ b)|	}.
\end{equation}
\end{definition}

\begin{lemma}\label{lem:number}
		Let  $B\subset A$ with $|A|-|B|\geq 2$; for $a, b \in A\setminus B$ the number of elements in set 
	$\mathrm{S}(a~ A\setminus B~b)$ is
	\begin{equation}\label{num}
	|\mathrm{S}(a~ A\setminus B~b) | = (|B|+1)! \times (|A\setminus B|-2)!
	\end{equation}
\end{lemma}
Indeed, we have $|B|!$ permutations for the elements in $B$ and there are $|B|+1$ ways to partition them into elements ranking above $a$ and below $b$; moreover, there are $|A\setminus B|-2)!$ permutations of the elements in $|A\setminus B|$ while keeping  $a$ and $b$ fixed. 

\subsection{Construction of  $\mathrm{P}[\,.\,]$ for $|A| = 4$}
We illustrate the construction  for  $A=\{a,b,c,d\}$. For  pairwise different and arbitrary $i, j, k, l \in A$, Table~\ref{tab1} contains all 12 possible rankings of the elements in $A$  where $i$ is ranked ``best''  and $j$ is ranked ``worst''  in $A\setminus B$.
\begin{table}[h]
	\begin{center}	
		$			\begin{array}{|c||c|c|c|c|} 
			\hline 	
			B& \emptyset & \{k\}& \{l\}   & \{k, l\}\\ 	\hline
			S(i ~A\setminus B~ j)	& S(i ~A~ j) & S(i~A\setminus \{k\}~j) &S(i~A\setminus \{l\}~j)  & S(i~A\setminus \{k, l\}~j) \\ 	\hline \hline
			& iklj &kilj  & likj &klij\\
			& ilkj & iljk& ikjl &  lkij\\
			&     &  &  & kijl\\
			\text{all rankings }	&     &  &  & lijk\\
			&     &  &  & ijkl\\
			&     &  &  & ijlk\\ \hline
			|S(i ~A\setminus B~ j)|	&  2!   & 2! & 2! &  3!\\\hline		
		\end{array}$\caption{Sets of rankings $S(i~A\setminus B~ j) \subset \Pi_A$, $i,j \in A\setminus B$. Note that interchanging the positions of $i$ and $j$ yields all 4! rankings .}
		\label{tab1}
	\end{center}
\end{table}	

Note that 
\begin{equation*}
	\sum_{B\in \mathscr{P}(A,\leq |A|-2)} (	|S(i ~A\setminus B~ j)|	+|S(j ~A\setminus B~ i)|	)= 4! =24 = |\Pi_A|.
 \end{equation*}

\begin{proposition}\label{pfor4}
	Let $A=\{i,j,k,l\}$ and  $B\in \mathscr{P}(A,\leq 2 )$; define for any $i,j\in A\setminus B$ ($i\neq j$)
	 \begin{equation*}\label{def}
		\mathrm{P}[\pi \in 	S(i~A\setminus B~j)] = \frac{K_{ij, B}}{|S(i~A\setminus B~ j)|	}.
	\end{equation*}
Then,
\begin{enumerate}
	\item[(a)] $\mathrm{P}(\Pi_A) = 1$;
	\item[(b)] 	$BW_{B}(i,j) = \mathrm{P}[S(i~B~j)] $.
\end{enumerate}
	
\end{proposition}
First, we have to show that $\mathrm{P}[\Pi_A]=1$. From Definition~\ref{def}, for each  $B\subset A$ 
\begin{align*}
	\sum_{\pi \in 	S(i ~A\setminus B~ j)} \mathrm{P}[\pi \in 	S(i ~A\setminus B~ j)] &= K_{ij, B}.
\end{align*}

We extend probability measure $\mathrm{P}(A)$ to a measure on subsets of $\Pi_A$ by setting
\[\mathrm{P}[S(i ~A\setminus B~ j)]= K_{ij,B}. \]
From Table~\ref{tab1} (or, Lemma~\ref{lemma:disjoint}), 

\[ \Pi_A =  \sum_{B\in \mathscr{P}(A,\leq |A|-2)} (	S(i ~A\setminus B~ j)	+S(j ~A\setminus B~ i)	),\]
so that 
\begin{align*}
	\mathrm{P}[\Pi_A]&= \sum_{B\in \mathscr{P}(A,\leq |A|-2)}  (K_{ij, B} +K_{ji, B}) \\
\end{align*}
\begin{align*}
\mathrm{P}[\Pi_A]&= \sum_{B\in \mathscr{P}(A,\leq |A|-2)}  K_{ij, B} +K_{ji, B} \\
&=(K_{ij,\emptyset}+K_{ji,\emptyset}) + (K_{ij,\{k\}}+K_{ji,\{k\}}) +	(K_{ij,\{l\}}+K_{ji,\{l\}})+	(K_{ij,\{k,l\}}+K_{ji,\{k,l\}}).
\end{align*}
Now, by definition
\begin{align}\label{emptyset}
&K_{ij,\emptyset}+K_{ji,\emptyset} = BW_A(i,j) + BW_A(j,i);\\
&K_{ij,\{k\}}+K_{ji,\{k\}}= BW_{A-\{k\}} (i,j) - BW_A(i,j) +BW_{A-\{k\}} (j,i) - BW_A(j,i);\\
&K_{ij,\{l\}}+K_{ji,\{l\}}= BW_{A-\{l\}} (i,j) - BW_A(i,j) +BW_{A-\{l\}} (j,i) - BW_A(j,i);\\
&	K_{ij,\{k,l\}}=BW_{\{i,j\}}(i,j)  - [ BW_{A-\{k\}} (i,j)+ BW_{A-\{l\}} (i,j)]   +BW_A(i,j);   \\
&K_{ji,\{k,l\}} = BW_{\{i,j\}}(j,i)  - [ BW_{A-\{k\}} (j,i)+ BW_{A-\{l\}} (j,i)]   +BW_A(j,i) .
\end{align}
 Abbreviating $BW_A(i,j) =(ij), BW_{A-\{k\}} (i,j)= (ij|k), BW_{A-\{l\}} (i,j)= (ij|l)$ and $BW_{\{i,j\}}(i,j)=(ij|k,l)$, collecting all terms from Equations  7--11,  reordering, and canceling yields

\begin{align*}
\mathrm{P}[\Pi_A]&=\\	
&(ij) + (ji) \\
&+ (ij|k)-(ij)+(ji|k)- (ji)\\
&+(ij|l)-(ij)+(ji|l)-(ji)\\
&+(ij|k,l)-[(ij|k)+(ij|l)]+(ij)\\
&+(ji|k,l)-[(ji|k)+(ji|l)]+(ji)\\
&=(ij) + (ji)-(ij)-(ji)-(ij)-(ji)+(ij) +(ji)\\
&  +(ij|k)+ (ji|k)-(ij|k)-(ji|k) +(ij|l)+ (ji|l)-(ij|l)-(ji|l)\\
&+(ij|k,l)+(ji|k,l)\\
&=1.
\end{align*}

Next, from Lemma~\ref{lem:BM}, we have to show that 
\begin{equation*}
	BW_{B}(a,b) = \mathrm{P}[S(a~B~b)] 
\end{equation*}
for all $B\in \mathscr{P}(A, \geq 2)$, 	and $a,b \in  B$. Equivalently, we show that 
\begin{equation*}
	BW_{A\setminus B}(a,b) = \mathrm{P}[S(a~A\setminus B~b)] 
\end{equation*}
for all $B\in \mathscr{P}(A, \le|A|- 2)$, 	and $a,b \in  A \setminus B$. 

Indeed, from Lemma~\ref{lemma:disjoint}  and  Theorem~\ref{moebius}
\begin{align*}
	\mathrm{P}[S(a~A\setminus B~b)]  &= \mathrm{P}\left[ \sum_{C\in \mathscr{P}( B)}\sum_{\pi\in \Pi_C} \sum_{(\pi_1 \pi_2)=\pi} S(\pi_1 a ~A\setminus C~ b \pi_2)\right]\\
	&=  \sum_{C\in \mathscr{P}( B)} \mathrm{P}\left[ \sum_{\pi\in \Pi_C} \sum_{(\pi_1 \pi_2)=\pi} S(\pi_1 a ~A\setminus C~ b \pi_2)\right]\\
	& = \sum_{C \in \mathscr{P}(B)}K_{ab, C}\\
	&= BW_{A\setminus B}(a,b).
\end{align*}
This completes the proof of Proposition~\ref{pfor4}.

\newpage
%

\bibliography{best-worst_corr}
\end{document}